
\input amstex

\magnification 1200
\loadmsbm
\parindent 0 cm

\input amssym
\input amssym.def

\define\nl{\bigskip\item{}}
\define\snl{\smallskip\item{}}
\define\inspr #1{\parindent=20pt\bigskip\bf\item{#1}}
\define\iinspr #1{\parindent=27pt\bigskip\bf\item{#1}}
\define\einspr{\parindent=0cm\bigskip}

\define\co{\Delta}

\define\ot{\otimes}

\define\tr{\triangleright}
\define\tl{\triangleleft}

%
%
\centerline{\bf Locally compact quantum groups.}
\centerline{\bf Radford's $S^4$ formula. \rm ($^+$)}
\bigskip
\centerline{\it A.\ Van Daele \rm ($^*$)}
\bigskip\bigskip
\bf Abstract \rm
\bigskip
Let $A$ be a finite-dimensional Hopf algebra. The left and the right integrals on $A$ are related by means of a distinguished group-like element $\delta$ of $A$. Similarly, there is this element $\widehat\delta$ in the dual Hopf algebra $\widehat A$. Radford showed that 
$$S^4(a)=\delta^{-1}(\widehat\delta\tr a \tl \widehat\delta^{-1})\delta$$
for all $a$ in $A$ where $S$ is the antipode of $A$ and where $\tr$ and $\tl$ are used to denote the standard left and right actions of $\widehat A$ on $A$. The formula still holds for multiplier Hopf algebras with integrals (algebraic quantum groups).
\smallskip
In the theory of locally compact quantum groups, an analytical form of Radford's formula can be proven (in terms of bounded operators on a Hilbert space).
\smallskip
In this talk, we do not have the intention to discuss Radford's formula as such, but rather to use it, together with related formulas, for illustrating various aspects of the road that takes us from the theory of Hopf algebras (including compact quantum groups) to multiplier Hopf algebras (including discrete quantum groups) and further to the more general theory of locally compact quantum groups.
\bigskip
\bigskip
\it August 2007 \rm (Version 1.0)
\bigskip
\bigskip
\vskip 5 cm
\hrule
\bigskip
\parindent 0.5 cm
\item{($^+$)} Notes for a talk at the {\it Workshop on Quantum Groups and Noncommutative Geometry} in Bonn (August 2007).
\item{($^*$)} Address: Department of Mathematics, K.U.\ Leuven, Celestijnenlaan 200B, B-3001 Heverlee, Belgium. E-mail: Alfons.VanDaele\@wis.kuleuven.be
\parindent 0 cm

\newpage


\bf 0. Introduction \rm
\nl
As we have mentioned in the abstract, this note is about different steps along the road from the (purely algebraic) theory of Hopf algebras to the (analytical) theory of locally compact quantum groups. The formula of Radford, under its different forms at each level, is only used to illustrate certain aspects in this development.
\nl
In {\it Section} 1, we start with the simplest case. We take a finite-dimensional Hopf algebra $A$ and we recall Radford's formula for the fourth power of the antipode in this case (see [R]), introducing the terminology that will be used further. We use $S$ for the antipode and $\delta$ and $\widehat \delta$ for the distinguished group-like elements in $A$ and the dual $\widehat A$. We call these the modular elements for reasons we explain later. We are also interested in the $^*$-algebra case and in particular when the underlying algebra is an operator algebra. This means that $A$ can be represented as a $^*$-algebra of operators on a (finite-dimensional) Hilbert space. Then however, the integrals are positive, the modular elements are $1$ and $S^2=\iota$ (the identity map) so that Radford's formula becomes a triviality. We speak about a finite quantum group but in the literature, it is usually called a finite-dimensional Kac algebra (see [E-S]).
\snl
In {\it Section} 2, we first consider the case of a Hopf algebra $A$, not necessarily finite-dimension\-al, but with integrals (a co-Frobenious Hopf algebra). Radford's formula in this case was obtained in [B-B-T] where the modular element $\widehat\delta$ is seen as a homomorphism from $A$ to $\Bbb C$. In this note however, we consider the dual $\widehat A$ of this Hopf algebra and describe it as a multiplier Hopf algebra. The element $\widehat \delta$ is then an element in the multiplier algebra $M(\widehat A)$ of $\widehat A$.
\snl
In the operator algebra framework, we get here (essentially) a compact quantum group (as introduced by Woronowicz in [W2] and [W3]). In this setting, we necessarily have $\delta=1$, but now it can happen that $\widehat \delta\neq 1$ (e.g.\ for the compact quantum group $SU_q(2)$, see [W1]). We also consider discrete quantum groups. They were first introduced in [P-W] as duals of compact quantum groups. Later they have been studied, as independent objects and independently in [E-R] and [VD2]. In this case of course, $\widehat\delta=1$ while possibly $\delta\neq 1$. Radford's formula gives $S^4(a)=\delta^{-1}a\delta$ for all $a$ in the algebra. In fact, one can define the square root $\delta^\frac12$ of $\delta$ and show that even $S^2(a)=\delta^{-\frac12}a\delta^\frac12$. It is a fundamental formula for discrete quantum groups.
\snl
{\it Section} 3 is about algebraic quantum groups. We already needed the notion of a multiplier Hopf algebra (see [VD1]) in Section 2 for properly dealing with discrete quantum groups. However, it is only in this section that we introduce the concept. We also look at the case with integrals and then we speak about algebraic quantum groups (cf.\ [VD3]). For an algebraic quantum group $(A,\co)$, it is possible to define a dual $(\widehat A,\widehat\co)$. It is again an algebraic quantum group. This duality extends the duality of finite-dimensional Hopf algebras (as used in Section 1), as well as the duality between compact and discrete quantum groups (as in Section 2). Also in this more general case, we have the existence of the modular elements $\delta$ and $\widehat\delta$, now in the multiplier algebras, and Radford's formula is still valid. It seems appropriate to give a proof (or rather sketch it) in this situation because it will follow easily from well-known results in the theory (see [D-VD-W]). As this case is more general than the previous ones (e.g.\ the finite-dimensional and the co-Frobenius Hopf algebras), this proof is also valid for these earlier cases.
\snl
Also here, we consider the $^*$-algebra case and in particular when the integrals are positive. Then, the underlying algebras are operator algebras (now $^*$-algebras of {\it bounded} operators on a possibly infinite-dimensional Hilbert space). We also have an analytical form of Radford's formula here and it is very close to the form we will obtain in the still more general case of locally compact quantum groups (in Section 4). Observe that now it can happen that both $\delta$ and $\widehat \delta$ are non-trivial.
\snl
It should not come as a surprise that, for $^*$-algebraic quantum groups, we can formulate a form of Radford's result that is similar to the one we will obtain for general locally compact quantum groups. After all, the theory of $^*$-algebraic quantum groups has been a source of inspiration for the development of locally compact quantum groups (as found in [K-V1], [K-V2] and [K-V3]). See e.g.\ the paper by Kustermans and myself [K-VD] and also the more recent paper entitled {\it Multiplier Hopf $^*$-algebras with positive integrals: A laboratory for locally compact quantum groups} [VD6].
\snl
Finally, in {\it Section} 4 we briefly discuss the most general and tecnically far more difficult case of a locally compact quantum group. We recall the definition (within the setting of von Neumann algebras) and we explain how the basic ingredients of the analytical form of Radford's result are constructed. About the proof, we have to be very short because this would take us too far. Nevertheless, we say something about it 
and especially, what kind of similarities there are with the case of algebraic quantum groups. Observe some differences in conventions in this section.
\nl
This note contains no new results. It is more like a short survey of various levels, from Hopf algebras to locally compact quantum groups, making a link between the purely algebraic approach to quantum groups and the operator algebra approach. It is well-known that working with operator algebras in this context puts sometimes very severe restrictions on possible results, special cases and examples. Think e.g.\ of the fact that it forces the square of the antipode to be the identity map in the finite-dimensional case (see Section 1). On the other hand, it also has some nice advantages like the analytic structure of a $^*$-algebraic quantum group (see Section 3). In any case, we are strongly convinced that a fair amount of knowlegde of 'the other side' can be of great help, not only for a basic understanding, but also because it sometimes provides different and handy tools to obtain new results or to treat old results in a better way. We think Radford's formula is a good illustration of this fact. Therefore, with this note, we hope to contribute to increase the interest of algebraists in the analytical aspects and vice versa.
\nl
Let us finish this introduction with some notation and conventions, as well as with providing some basic references. More of this will be given throughout the note.
\snl
We work with associative algebras over the comlex numbers since we often will also consider an involution on the algebra, making it into an operator algebra. The algebras need not have an identity, but we always assume that the product, as a bilinear map, is non-degenerate. This allows to consider the algebra as a two-sided ideal sitting in the multiplier algebra. If the algebra has a unit, we denote it by $1$. This will also be used for the unit in the multiplier algebra. We will systematically use $\iota$ for the identity map. 
\snl
We use $A'$ for the space of all linear functionals on a vector space $A$ and call it the dual space of $A$. Often, we will consider a suitable subspace of this full dual space. Most of the time, our tensor products are purely algebraic, except in the last section on locally compact quantum groups where we work with von Neumann algebras and von Neumann algebraic tensor products. Unfortunately, some other conventions in Section 4 are also different from those in the earlier sections. This is mainly due to differences between the algebraic approach and the operator algebra approach.
\snl
The basic references for Hopf algebras are of course [A] and [S]. For compact quantum groups we have [W2] and [W3], see also [M-VD]. For discrete quantum groups we refer to [P-W], [E-R] and [VD2]. The basic theory of multiplier Hopf algebras is found in [VD1] and when they have integrals, the reference is [VD3]. See also [VD-Z] for a survey paper on the subject. Finally, the general theory of locally compact quantum groups is developed in [K-V1], [K-V2] and [K-V3]. See also [M-N] and [M-N-W] for a different approach and [VD8] for a more recent and simpler treatment of the theory.
\nl
\nl
\bf Acknowledgements \rm
I would like to thank the organizers of the Workshop on Quantum Groups and Noncommutative Geometry (MPIP Bonn, August 2007) for giving me the opportunity to talk about this subject. I am also grateful to P.M.\ Hajac who drew my attention to the paper by Kaufman and Radford [K-R].
\nl
\nl
%
%
\bf 1. Finite quantum groups \rm
\nl
Let $A$ be a finite-dimensional Hopf algebra (over the complex numbers) with coproduct $\co$, counit $\varepsilon$ and antipode $S$. Let $\widehat A$ denote the dual Hopf algebra of $A$. We will us the {\it pairing notation}. So, if $a\in A$ and $b\in \widehat A$ we write 
$\langle a,b \rangle$ for the value of $b$ in the element $a$.
\snl
Let $\varphi$ be a left integral on $A$. There exists a distinguished group-like element $\delta$ in $A$ defined by the formula $\varphi(S(a))=\varphi(a\delta)$ for all $a\in A$. We will call $\delta$ the {\it modular element} (for reasons we will explain later, in Section 3). Similarly, when $\widehat\varphi$ is a left integral on $\widehat A$, there is the modular element $\widehat \delta$ in $\widehat A$ satisfying $\widehat\varphi(S(b))=\widehat\varphi(b\widehat\delta)$ for all $b\in \widehat A$. 
\snl
Now we can state Radford's formula (see [R]):

\inspr{1.1} Theorem \rm
For all $a\in A$, we have  
$$S^4(a)=\delta^{-1}(\widehat\delta\tr a \tl \widehat\delta^{-1})\delta.$$
\einspr

We use the standard left and right actions of the dual $\widehat A$  on $A$ defined by
$$b\tr a= \sum_{(a)} a_{(1)}\langle a_{(2)},b\rangle \quad\qquad \text{and} \quad\qquad
a\tl b= \sum_{(a)} a_{(2)}\langle a_{(1)},b\rangle$$
for $a\in A$ and $b\in \widehat A$ (where we use the Sweedler notation).
\einspr

Later, we will give a proof of this formula in the more general setting of algebraic quantum groups (see Section 3).
\nl
Let us also consider the case of a Hopf $^*$-algebra. We assume that $A$ is a $^*$-algebra and that $\co$ is a $^*$-homomorphism. Then $\varepsilon$ is a $^*$-homomorphism but $S$ need not be a $^*$-map. In stead, it is invertible and satisfies $S(a)^*=S^{-1}(a^*)$ for all $a$. So, it is a $^*$-map if and only if $S^2=\iota$, the identiy map.
\snl
If moreover $A$ is an operator algebra, then there exists a {\it positive} left integral $\varphi$ (and conversely). Then necessarily $\varphi(1)>0$ so that left and right integrals coincide. This implies that $\delta=1$. One can show that again $\widehat A$ will be an operator algebra and so also $\widehat\delta=1$. Radford's formula implies that in this case $S^4=\iota$. In fact, it follows that already $S^2=\iota$ and that the integrals are traces. We will give a short argument later in the more general case of a discrete quantum group (see the next section and also Section 3).
\snl
In this note, we will call a finite-dimensional Hopf $^*$-algebra with positive integrals a {\it finite quantum group}. In the literature however, one often calls it a finite-dimensional Kac algebra (see [E-S]).
\nl\nl
%
%
\bf 2. Compact and discrete quantum groups \rm
\nl
Now, let $A$ be any Hopf algebra. We do no longer assume that it is finite-dimensional, but we require that it has integrals. Assume also that it has an invertible antipode. Again there exists a unique group-like element $\delta$ in $A$ such that $\varphi(S(a))=\varphi(a\delta)$ for all $a\in A$ when $\varphi$ is a left integral on $A$.
\snl
The dual space $A'$ is an algebra but no longer a Hopf algebra (in general). However, there still is the distinguished element $\widehat\delta\in A'$. It is a homomorphism, it is invertible and Radford's formula is still valid.
For all $a\in A$, we have  
$$S^4(a)=\delta^{-1}(\widehat\delta\tr a \tl \widehat\delta^{-1})\delta.$$
The actions are defined as before by
$$f\tr a= \sum_{(a)} f(a_{(2)}) a_{(1)} \quad\qquad \text{and} \quad\qquad
a\tl f= \sum_{(a)} f(a_{(1)})a_{(2)}$$
for all $a\in A$ and $f\in A'$.
\snl
The proof we plan to give later (for algebraic quantum groups) will also include this case.
\nl
If moreover $A$ is a $^*$-algebra and $\co$ a $^*$-homomorphism, still $\varepsilon$ will be a $^*$-homomorphism and $S(a)^*=S^{-1}(a^*)$ for all $a\in A$. And if $A$ is an operator algebra, the left integral is positive, it is also a right integral and so $\delta=1$.
\snl
We agree to use the term {\it compact quantum group} for this case. Indeed, it is essentially a compact quantum group as defined by Woronowicz in [W3].
\snl
Remark that $\widehat \delta$ need not be $1$ in this case, the integrals need not be traces and $S^2\neq \iota$ is still possible. The standard example where this happens is the quantum $SU_q(2)$ (see [W1]).
\nl
Let us now consider the case of a discrete quantum group. Discrete quantum groups can be obtained as duals of compact quantum groups. Although it is more natural to treat them within the framework of multiplier Hopf algebras (see later), we will briefly consider the case already now (and see why we need to pass to multiplier Hopf algebras).
\snl
The following result is part of the motivation for what we will do later.

\inspr{2.1} Proposition \rm
Let $A$ be a Hopf algebra with a left integral $\varphi$. Define the dual $\widehat A$ as the subspace of $A'$ containing all elements of the form $\varphi(\,\cdot\,a)$ with $a\in A$. It is a subalgebra of $A'$. If we define the coproduct $\widehat\co:A'\to (A\ot A)'$ by dualizing the product on $A$, we find that
$$\widehat\co(\widehat A)(1\ot \widehat A)\subseteq \widehat A\ot \widehat A \quad\qquad\text{and}\quad\qquad
(\widehat A\ot 1)\widehat\co(\widehat A)\subseteq \widehat A\ot \widehat A$$
in the algebra $(A\ot A)'$. 
\einspr

So, we get that $\widehat\co$ maps $\widehat A$ into the multiplier algebra $M(\widehat A\ot \widehat A)$ (as we will define it later). Moreover, the pair $(\widehat A,\widehat\co)$ is a multiplier Hopf algebra (and not a Hopf algebra in general).
\snl
If we define $\widehat\psi(b)=\varepsilon(a)$ when $b=\varphi(\,\cdot\,a)$, we get a right integral on $\widehat A$. This means here that
$$(\widehat\psi\ot\iota)(\widehat\co(b)(1\ot b'))=\widehat\psi(b)b'$$
for all $b,b'\in \widehat A$. The antipode $S$ leaves $\widehat A$ invariant and converts $\widehat\psi$ to a left integral $\widehat\varphi$ on $\widehat A$. The element $\widehat \delta$, considered earlier, is in $M(\widehat A)$ and still satifsies $\widehat\varphi(S(b))=\widehat\varphi(b\widehat\delta)$ for all $b\in \widehat A$.
\nl
If $A$ is a compact quantum group, it turns out that $\widehat A$ is a direct sum of matrix algebras. This takes us to the following definition of a discrete quantum group.

\inspr{2.2} Definition \rm
A {\it discrete quantum group} is a pair $(A,\co)$ where $A$ is a direct sum of matrix algebras (with the standard involution), $\co$ is a coproduct on $A$ and such that there is a counit $\varepsilon$ and an antipode $S$.
\einspr

It is not a Hopf algebra (except when it is a finite direct sum), but it is a multiplier Hopf algebra (see further). Indeed, we have $\co(A)\subseteq M(A\ot A)$, the multiplier algebra of $A\ot A$, but in general $\co(A)$ does not belong to $A\ot A$ itself.
\nl
For discrete quantum groups, we can prove (among other things) the following result.

\inspr{2.3} Theorem \rm
There exists a positive left integral $\varphi$ and a positive group-like element $\delta$ in the multiplier algebra $M(A)$ of $A$ defined by $\varphi(S(a))=\varphi(a\delta)$ for all $a\in A$. This element moreover satisfies
$$S^2(a)=\delta^{-\frac12}a\delta^{\frac12}$$
for all $a$. We also have $\varphi(ab)=\varphi(bS^2(a))$ for all $a,b\in A$ and therefore, the map $a\mapsto \varphi(a\delta^\frac12)$ is a trace on $A$.
\einspr

The first formula is a slightly stronger version of Radford's formula for these discrete quantum groups. It can be dualized to get a similar expression for the square $S^2$ of the antipode of a compact quantum group. 
\snl
One way to develop discrete quantum groups is by viewing them as duals of compact quantum groups (as done in [P-W]). This however is not the best choice. It is relatively easy to develop the theory of discrete quantum groups (and prove the above results) directly from the definition above. Using the standard trace on each component, one can obtain quickly a formula for both integrals as well as for the modular element. See e.g.\ [VD2].
\snl
It can happen that $\delta\neq 1$ (so that left and right integrals are different). It can also happen that $S^2\neq \iota$ so that the integrals are not traces. This can of course only happen if $\delta\neq 1$.
\snl
The standard example is the dual of the compact quantum group $SU_q(2)$ whose underlying algebra is the direct sum $\oplus_{n=0}^\infty M_n$. All objects can easily be given in terms of the deformation parameter $q$, except for the comultiplication (which is quite complicated), see e.g.\ [VD4].
\snl
On the other hand, if $\delta=1$ we must have that $S^2=\iota$ and that the integrals are traces. This generalizes the corresponding result for finite quantum groups as we have seen in Section 1. Observe also that if we have a quantum group that is both discrete and compact, it must be a finite quantum group.
\nl
\nl
%
%
\bf 3. Algebraic quantum groups \rm
\nl
Discrete and compact quantum groups are special cases of algebraic quantum groups. Also the duality of algebraic quantum groups generalizes the one between discrete and compact quantum groups. We will briefly review this theory. For details, we refer to the literature, see [VD1], [VD3] and [VD-Z].
\snl
The basic ingredient is that of a multiplier Hopf algebra:

\inspr{3.1} Definition \rm
Let $A$ be an algebra over $\Bbb C$, with or without identity, but with a non-degenerate product. A {\it coproduct} (or comultiplication) on $A$ is a non-degenerate homomorphism $\co :A \to M(A\ot A)$ (the multiplier algebra of $A\ot A$), satisfying coassociativity $(\co\ot\iota)\co=(\iota\ot\co)\co$. The pair $(A,\co)$ is called a {\it (regular) multiplier Hopf algebra} if there exists a counit and an (invertible) antipode. If $A$ is a $^*$-algebra and $\co$ a $^*$-homomorphism, regularity is automatic and we call it a multiplier Hopf $^*$-algebra.  
\einspr

There is a lot to say about this definition and we refer to the literature for details. However, it is important to notice that any Hopf ($^*$-)algebra is a multiplier Hopf ($^*$-)algebra and conversely, if the underlying algebra of a multiplier Hopf algebra has an identity, it is actually a Hopf algebra. Also remark that the counit and the antipode are unique.
\snl
Next, we consider algebraic quantum groups:

\inspr{3.2} Definition \rm 
Let $(A,\co)$ be a regular multiplier Hopf algebra. A {\it left integral} is a non-zero linear functional $\varphi:A\to \Bbb C$ satisfying left invariance $(\iota\ot\varphi)\co(a)=\varphi(a)1$ in $M(A)$ for all $a\in A$. Similarly, a right integral is defined. 
\einspr

If a left integral $\varphi$ exists, also a right integral $\psi$ exists (namely $\psi=\varphi\circ S$). In that case, we call $(A,\co)$  an {\it algebraic quantum group}. If moreover $(A,\co)$ is a multiplier Hopf $^*$-algebra with a positive left integral $\varphi$ (i.e.\ such that $\varphi(a^*a)\geq 0$ for all $a$), then also a positive right integral exists (which is not a trivial result!). In that case, we call $(A,\co)$ a $^*$-algebraic quantum group.
\snl
Remark that the term 'algebraic' does not refer to the possible quantization of algebraic groups, but we use it rather because $^*$-algebraic quantum groups are locally compact quantum groups (considered in the next section) that can be treated with purely algebraic techniques. 
\snl
Integrals on regular multiplier Hopf algebras are unique (up to a scalar) if they exist. They are faithful in the sense that (for the left integral $\varphi$) we have $a=0$ if either $\varphi(ab)=0$ for all $b$ or $\varphi(ba)=0$ for all $b$. From the uniqueness it follows that there is a constant $\nu$ (the {\it scaling constant}), given by $\varphi(S^2(a))=\nu \varphi(a)$ for all $a\in A$. It can happen that $\nu\neq 1$ but when $A$ is a $^*$-algebraic quantum group (with positive integrals), we must have $\nu=1$ (see [DC-VD]).
\snl
In general, integrals need not be traces, but there exist automorphisms $\sigma$ and $\sigma'$ (called the {\it modular automorphisms}) satisfying 
$$\varphi(ab)=\varphi(b\sigma(a))\qquad \qquad \qquad \psi(ab)=\psi(b\sigma'(a))$$
for all $a,b\in A$ when $\varphi$ is a left integral and $\psi$  a right integral. The term 'modular' comes from operator algebra theory and the modular automorphism group of a faithful normal state (or semi-finite weight) on a von Neumann algebra (see the next section).
\snl
Important for us in this note that focuses on Radford's formula is the {\it modular element} $\delta$. It is a group-like element in the multiplier algebra $M(A)$ satisfying $\varphi(S(a))=\varphi(a\delta)$ for all $a$ just as in the case of Hopf algebras with integrals. It can be defined, using the uniqueness of integrals, by the formula $(\varphi\ot\iota)\co(a)=\varphi(a)\delta$ for all $a$. In this case, the term 'modular' is used because it is related with the modular function for a non-unimodular locally compact group. In fact, also the modular automorphism group in the theory of von Neumann algebras finds its origin in the theory of non-unimodular locally compact groups.
\snl
There are many relations among these objects and again, we refer to the literature. 
\nl
For any algebraic quantum group, we have a dual:

\inspr{3.3} Theorem \rm
Let $(A,\co)$ be an algebraic quantum group. Define the subspace $\widehat A$ of the dual space $A'$ of functionals of the form $\varphi(\,\cdot\,a)$ where $a\in A$. The adjoints of the coproduct and the product on $A$ define a product and a coproduct $\widehat\co$ on $\widehat A$, making $(\widehat A,\widehat\co)$ into an algebraic quantum group, called the {\it dual} of $(A,\co)$. A right integral $\widehat\psi$ on $\widehat A$ is given by the formula $\widehat\psi(\omega)=\varepsilon(a)$ when $\omega=\varphi(\,\cdot\,a)$ and $a\in A$. If $(A,\co)$ is a $^*$-algebraic quantum group, then so is $(\widehat A,\widehat\co)$ and $\widehat\psi$ as defined before is positive when $\varphi$ is positive. 
\einspr

The last statement in the above theorem is a consequence of {\it Plancherel's formula}. Here it says that $\widehat\psi({\widehat a}^*\widehat a)=\varphi(a^*a)$ if $a\in A$ and $\widehat a=\varphi(\,\cdot\, a)$, its {\it Fourier transform}. 
\snl
Also remakr that the dual of $(\widehat A,\widehat\co)$ is again $(A,\co)$.   
\snl
We will use the pairing notation (as we have already done in Section 1 for a finite-dimensional Hopf algebra and its dual). We also have the standard actions of $\widehat A$ on $A$ and of $A$ on $\widehat A$. In the first case, we have 
$$\align \langle b\tr a,b'\rangle &= \langle a, b'b\rangle \\
         \langle a\tl b,b'\rangle &= \langle a, bb'\rangle
\endalign$$
for all $a\in A$ and $b,b'\in B$. It is not completely obvious that these elements are well-defined in $A$, but it can be shown. Moreover, these actions are unital. This means that elements of the form $b\tr a$ with $a$ in $A$ and $b$ in $\widehat A$ span all of $A$ and similarly for the right action. See [Dr-VD] and [Dr-VD-Z].
\nl
Let us now first state some of the formulas relating the various objects of $(A,\co)$ and indicate how they can be proven. We use the notations introduced before.

\inspr{3.4} Proposition \rm
Let $(A,\co)$ be an algebraic quantum group. We have $\sigma\circ S\circ\sigma'=S$ and $\delta \sigma(a)=\sigma'(a)\delta$ and for all $a$. Also for all $a\in A$
we have
$$\co(\sigma(a))=(S^2\ot \sigma)\co(a).$$
\einspr

The first formulas follow in a straighforward way from the definitions of $\sigma$ and $\sigma'$ with $\psi=\varphi\circ S=\varphi(\,\cdot\,\delta)$. For the second one, we use that for all $a,b\in A$,
$$S((\iota\ot\varphi)(\co(a)(1\ot b)))=(\iota\ot\varphi)((1\ot a)\co(b)),$$
two times in combination with the definition of $\sigma$. This last formula itself follows easily from left invariance of $\varphi$ and the standard properties of the antipode.
\snl
We will also need some other properties. We have that the automorphisms $S^2$, $\sigma$ and $\sigma'$ all commute with each other. And we also have that $\sigma(\delta)=\sigma'(\delta)=\frac{1}{\nu}\delta$ where $\nu$ is the scaling constant.
\nl
Next, we state and prove some of the formulas relating objects of $(A,\co)$ with objects of the dual $(\widehat A,\widehat\co)$.

\inspr{3.5} Proposition \rm
Let $(A,\co)$ be an algebraic quantum group and let $(\widehat A,\widehat\co)$ be its dual. We have $\widehat\delta^{-1}=\varepsilon\circ\sigma$ where $\widehat\delta$ is the modular element of $\widehat A$, seen as a homomorphism of $A$. Also $\sigma(a)=\widehat\delta^{-1}\tr S^2(a)$ for all $a\in A$.
\snl\bf Proof\rm:
To prove the first formula, we start with $c\in A$ and we take the element $b=\varphi(\,\cdot\, c)$ in the dual $\widehat A$. Because for all $a,a'$ in $A$ we have $\varphi(a'c\sigma(a))=\varphi(aa'c)$, we get $\varphi(\,\cdot\, c\sigma(a))=b\tl a$. If we apply $\widehat \psi$ we find $\varepsilon(c\sigma(a))=\widehat\psi(b\tl a)$. Because $(\iota\ot\widehat\psi)\widehat\co(b)=\widehat\psi(b)\widehat\delta^{-1}$ (a formula that can easily be obtained from the definition of $\widehat\delta$  by using the antipode), we get $\widehat\psi(b\tl a)=\widehat\psi(b)\langle a,\widehat\delta^{-1}\rangle$. Combining all results and using that $\widehat\psi(b)=\varepsilon(c)$, we find the first formula of the proposition.
\snl
To obtain the second formula, consider the equation $\co(\sigma(a))=(S^2\ot \sigma)\co(a)$, obtained in the previous proposition, apply $\iota\ot\varepsilon$  and use the first formula of this proposition.
\einspr 

In the proof above, we have used the left action of $A$ on $\widehat A$. We also have looked at $\widehat\delta^{-1}$ as a linear functional on $A$ by extending the pairing between $A$ and $\widehat A$ to $M(\widehat A)$ in an obvious way. If the quantum group is counimodular, that is if $\widehat\delta=1$, it follows from these results that $\sigma=S^2$. This is the case for discrete quantum groups as we saw in Theorem 2.3
\nl
Now we are ready to give a simple proof of Radford's formula for algebraic quantum groups.

\inspr{3.6} Theorem \rm
Let $(A,\co)$ be an algebraic quantum group. When $\delta$ and $\widehat\delta$ are the modular elements in $A$ and its dual $\widehat A$, then 
$$S^4(a)=\delta^{-1}(\widehat\delta\tr a \tl \widehat\delta^{-1})\delta$$
for all $a\in A$.
\snl \bf
Proof\rm:
From the second formula in Proposition 3.5 we find $\widehat\delta\tr a=S^2(\sigma^{-1}(a))$. Similarly, or by applying the antipode on this formula, we obtain $a \tl \widehat\delta^{-1} = S^2({\sigma'}(a))$. If we combine these two formulas with the relation $\sigma'(a)=\delta\sigma(a)\delta^{-1}$ and use that $S^2(\delta)=\delta$, we get Radford's formula.
\einspr 

The proof we have given can be found in [D-VD-Z] and in [D-VD], where we have generalized this result further to algebraic quantum hypergroups.

\nl
Next, let us look at the case of a $^*$-algebraic quantum group. The requirement of positivity of the integrals is quite strong. We have mentioned already that it forces the scaling constant $\nu$ to be 1. On the other hand, we end up with an operator algebra and this allows to work on Hilbert spaces and use spectral theory. In this case, we arrive at what is called the {\it analytic structure} of a $^*$-algebraic quantum group (see [K] and also [DC-VD]). Roughly speaking, it means that powers of $S^2$, $\sigma$, $\sigma'$ and $\delta$ all have analytical extensions to the whole complex plane. More precisely, we get the following result. We only consider $S^2$ and $\delta$ because we focus in this note on Radford's formula.

\inspr{3.7} Proposition \rm
Let $(A,\co)$ be a $^*$-algebraic quantum group. There exists an analytic function $\tau:z\mapsto \tau_z$ on $\Bbb C$ such that $\tau_z$ is an automorphism of $(A,\Delta)$, that $\tau_{z+y}=\tau_z\circ\tau_{y}$ for all $z,y\in \Bbb C$ and so that $S^2=\tau_{-i}$. Similarly, there is an analytic function $z \mapsto \delta^z$ so that $\delta^z\in M(A)$, that $\delta^{z+y}=\delta^z \delta^y$ for all $z,y\in \Bbb C$ and such that $\delta^z=\delta$ for $z=1$ (justifying the notation).
\einspr

Analyticity here is in a strong sense. In the first case, we want $z\mapsto f(\tau_z(a))$ analytic for all $a\in A$  and all $f\in A'$. In the second case, we want e.g.\ $z\mapsto f(a\delta^z)$ analytic for all $a\in A$ and $f\in A'$. These analytical extensions are unique.
\snl
Then, we can get the {\it analytical form} of Radford's formula. For real numbers, we obtain the following:

\inspr{3.8} Theorem \rm
Let $(A,\co)$ be a $^*$-algebraic quantum group. Let $\tau_z$ and $\delta^z$ for $z\in \Bbb C$ be defined as in the previous proposition. Consider also $\widehat\delta^z\in M(\widehat A)$ in a similar way. Then, for all $t\in \Bbb R$, we have
$$\tau_{2t}(a)=\delta^{-it}(\widehat\delta^{it}\tr a \tl \widehat\delta^{-it})\delta^{it}$$
for all $a\in A$.
\einspr

This is the form of Radford's formula that we will be able to generalize to general locally compact quantum groups (see the next section). The result however is true for all complex numbers. In particular, we can take $z=-\frac{i}{2}$. This yields
$$S^2(a)=\delta^{-\frac12}(\widehat\delta^\frac12\tr a \tl \widehat\delta^{-\frac12})\delta^\frac12$$
for all $a\in A$. Indeed, as a consequence of the result in Proposition 3.7, we can also define the square roots $\delta^\frac12$ and $\widehat\delta^\frac12$ in $M(A)$ and $M(\widehat A)$ respectively. These are still group-like elements.
\snl
We should make a reference to a paper by Kaufman and Radford here [K-R]. They discover the formula with the square roots for Drinfel'd doubles that are ribbon Hopf algebras.
\snl
Finally, consider some special cases. If e.g.\ $(A,\co)$ is counimodular, this is by definition when left and right integrals on $\widehat A$ are the same, so that $\widehat\delta=1$, we find that $S^2(a)=\delta^{-\frac12}a\delta^\frac12$ for all $a$. This is the formula that we have seen in Theorem 2.3 for discrete quantum groups. They are counimodular because compact quantum groups are unimodular. If $(A,\co)$ is both unimodular and counimodular, then we must have $S^2=\iota$. In this case, it follows from Proposition 3.5 that both $\sigma$ and $\sigma'$ are trivial. This means that the integrals are traces. This, in particular, applies to the case of finite quantum groups (as in Section 1).
\nl\nl
%
%
\bf 4. Locally compact quantum groups \rm
\nl
We start this section with the definition of a locally compact quantum group in the von Neumann algebra setting.

\inspr{4.1} Definition \rm
Let $M$ be a von Neumann algebra. A {\it coproduct} on $M$ is a normal unital $^*$-homomorphism $\co:M\to M\ot M$, the von Neumann algebraic tensor product, satisfying coassociativity $(\co \ot \iota)\co=(\iota\ot\co)\co$. If there exist faithful normal semi-finite weights $\varphi$ and $\psi$ on $M$ that are left, resp.\ right invariant, then the pair $(M,\co)$ is called a {\it locally compact quantum group}.
\einspr

We collect some important remarks about this concept:

\inspr{4.2} Remarks \rm 
i) The adapted form of continuity of $\co$ in the von Neumann algebra setting is expressed in the requirement that the coproduct is normal. \newline
ii) By this continuity, the $^*$-homomorphisms $\co\ot\iota$ and $\iota\ot\co$ are well-defined from $M\ot M$ to $M\ot M\ot M$ and so coassociativity makes sense.\newline
iii) A weight on a von Neumann algebra is, roughly speaking, an unbounded positive linear functional. It is called semi-finite if it is bounded on enough elements. And again it is called  normal if it satisfies the proper continuity.\newline
\einspr

For the theory of von Neumann algebras and the notions needed above, we refer to the books of Takesaki [T1] and [T2].
\snl
The weight $\varphi$ is called left invariant if $\varphi((\omega\ot\iota)\co(x))=\varphi(x)\omega(1)$ whenever $x$ is a positive element in the von Neumann algebra with $\varphi(x)<\infty$ and $\omega$ is a positive element in the predual $M_*$ of $M$. Similarly, right invariance of the weight $\psi$ is defined. These weights are unique (up to a scalar) and are called the left and right {\it Haar weights}. They are of course the analogues of the left and right integrals in the theory of $^*$-algebraic quantum groups.
\snl
It is also possible to define locally compact quantum groups in the framework of C$^*$-algebras, but that is somewhat more complicated. In fact, both approaches are equivalent in the sense that they define the same objects. We refer to the original works by Kustermans and Vaes; see [K-V1], [K-V2] and [K-V3]. Independently, the notion was also developed by Masuda, Nakagami and Woronowicz; see [M-N] and [M-N-W]. A more recent and simpler development of the theory can be found in [VD8] and a discussion on the equivalence of the C$^*$-approach and the von Neumann approach is e.g.\ given in [VD7].
\snl
The basic examples come from a locally compact group $G$. On the one hand, there is the abelian von Neumann algebra $L^\infty(G)$, defined with respect to the left Haar measure. The coproduct $\co$ is given as before by $\co(f)(p,q)=f(pq)$ when $f\in L^\infty(G)$ and $p,q\in G$. The invariant weights $\varphi$ and $\psi$ are obtained by integration with respect to the left and right Haar measures on the group. On the other hand, there is the group von Neuman algebra $VN(G)$ generated by the left regular representation $\lambda$ of the group on the Hilbert space $L^2(G)$. In this case, the coproduct is given by $\co(\lambda_p)=\lambda_p\ot \lambda_p$. The left and right integrals are the same. Formally, we must have $\varphi(\lambda_p)=0$, except when $p=e$, the identity of the group, but it is not so easy to define this weight properly.
\nl
Any multiplier Hopf $^*$-algebra with positive integrals, i.e.\ a $^*$-algebraic quantum group, gives rise to a locally compact quantum group (see [K-VD]):

\inspr{4.3} Theorem \rm
Let $(A,\co)$ be a $^*$-algebraic quantum group with left integral $\varphi$. Consider the GNS-representation $\pi_\varphi$ of $A$ associated with $\varphi$. The coproduct on $A$ yields a coproduct on the von Neumann algebra $M$ generated by $\pi_\varphi(A)$ making it into a locally compact quantum group.     
\einspr

The Haar weights are of course nothing else but the unique normal extensions of the original left and right integrals.
\snl
It is an interesting, but open problem to describe those locally compact quantum groups that can arise from $^*$-algebraic quantum groups as above. In the case of locally compact groups, the problem has been solved in [L-VD]. The requirement is that there exists a compact open subgroup. In particular, when $G$ is a totally disconected locally compact group, the two associated locally compact quantum groups are essentially $^*$-algebraic quantum groups. In connection with this problem, let us also observe the following. For any $^*$-algebraic quantum group, the scaling constant $\nu$ is necessarily $1$ (see [DC-VD]). However, there are examples of locally compact quantum groups where this is not the case (see [VD5]). We will come back to this statement later.
\nl
Let us now indicate how the theory of locally compact quantum groups is developed (as e.g.\ in [VD8]) and focus on the relevant formulas, needed to formulate Radford's result.
\snl
So, we start with a locally compact quantum group $(M,\co)$ with left and right Haar weights $\varphi$ and $\psi$ as in Definition 4.1. We recall the GNS construction:

\inspr{4.4} Proposition \rm
Denote by $\Cal N_\varphi$ the set of elements $x\in M$ so that $\varphi(x^*x)<\infty$. It is a dense left ideal of $M$ and $\varphi$ has a unique extension (still denoted by $\varphi$) to the $^*$-algebra spanned by elements of the form $x^*y$ with $x,y\in \Cal N_\varphi$. There exists a Hilbert space $\Cal H_\varphi$ and an injective linear map $\Lambda_\varphi:\Cal N_\varphi\to \Cal H_\varphi$ with dense range such that $\langle\Lambda_\varphi(x),\Lambda_\varphi(y)\rangle=\varphi(y^*x)$ for all $x,y\in \Cal N_\varphi$. There also exists a faithful, unital and normal $^*$-representation $\pi_\varphi$ of $M$ on $\Cal H_\varphi$ given by $\pi_\varphi(y)\Lambda_\varphi(x)=\Lambda_\varphi(yx)$ whenever $x\in \Cal N_\varphi$ and $y\in M$.
\einspr

In what follows, we will drop the index $\varphi$ and use $\Cal H$ and $\Lambda$ in stead of $\Cal H_\varphi$ and $\Lambda_\varphi$. We will also omit $\pi_\varphi$ and assume that $M$ acts directly on the space $\Cal H$.
\snl
Next, we recall some results from the Tomita-Takesaki modular theory (see e.g.\ [T2]):

\inspr{4.5} Proposition \rm
There is a closed, conjugate linear, possibly unbounded but densely defined involutive operator $T$ on $\Cal H$ so that $\Lambda(x)\in \Cal D(T)$, the domain of $T$, for any $x \in \Cal N_\varphi \cap \Cal N_\varphi^*$ and $T\Lambda(x)=\Lambda(x^*)$. If $T=J\nabla^\frac12$ denotes the polar decomposition of $T$, then $J$ is a conjugate linear isometric involutive operator and $\nabla$ a positive non-singular self-adjoint operator. If $M'$ denotes the commutant of $M$, we have $JMJ=M'$. Also $\nabla^{it}M\nabla^{-it}=M$ for all $t\in \Bbb R$.  
\einspr

It follows from this result that we can define a one-parameter group $(\sigma_t)$ of automorphisms of $M$, called the {\it modular automorphism group}, by $\sigma_t(x)=\nabla^{it}x\nabla^{-it}$ for $x\in M$ and $t\in \Bbb R$. A similar construction will give the modular automorphisms $(\sigma'_t)$ associated with the right Haar weight $\psi$.
\snl
Using a proper notion of an analytic extension, one can show that $\varphi(xy)=\varphi(y\sigma_{-i}(x))$ for the appropriate elements $x$ and $y$. So $(\sigma_{-i})$ plays the role of the modular automorphism $\sigma$ as we have it for $^*$-algebraic quantum groups. Similarly $\sigma'_{-i}$ plays the role of the modular automorphism $\sigma'$. We apologize for the possible confusion caused by the difference in notations used here (and further in this section). 
\snl
There is also something called the 'relative modular theory' when two weights are considered. If we apply results from this theory to the invariant weights $\varphi$ and $\psi$, we find the following:

\inspr{4.6} Proposition \rm
There exists a positive non-singular self-adjoint operator $\delta$ on the Hilbert space $\Cal H$ such that for all $t\in\Bbb R$ we have $\delta^{it}\in M$ and $\psi=\varphi(\delta^\frac12 \,\cdot\,\delta^\frac12)$.
\einspr

It should be mentioned that it is not so easy to interprete this last formula in a correct way.
\snl
When thinking of a $^*$-algebraic quantum group, where we have $\sigma(\delta)=\delta$  (because the scaling constant is trivial), we see that this formula is another form of the one we have for algebraic quantum groups, namely $\psi=\varphi(\,\cdot\,\delta)$. Here, we call $\delta$ the {\it modular operator}.
\nl
These are the first main ingredients of the theory. Remark that these objects are only dependent on the weights $\varphi$ and $\psi$ on the von Neumann algebras $M$ and seem in no way related with the coproduct structure. This is not completely correct as the result in Proposition 4.6 would not be true for any pair of weights.
\nl
Next, let us consider the dual locally compact quantum group $(\widehat M,\widehat \co)$ with left and right Haar weights $\widehat \varphi$ and $\widehat \psi$. The precise construction is quite involved but in essence, it is a careful analytic version of the same construction for $^*$-algebraic quantum groups.
\snl
The Hilbert space associated with the dual left Haar weight $\widehat\varphi$ is identified with $\Cal H$ and the map $\widehat\Lambda$ associated with $\widehat\varphi$ is defined in such a way that $\widehat\Lambda(\widehat x)=\Lambda(x)$ when $x$ is an appropriate element in $M$ and $\widehat x$ its Fourier transform $\varphi(\,\cdot\,x)$. Remark that a different convention is used in the sense that the dual coproduct is flipped causing, among other things, that the dual right integral $\widehat\psi$ is now the dual left integral $\widehat\varphi$. This convention is common in the operator algebra approach.
\snl
And just as for the original locally compact quantum group $(M,\co)$, we also have the conjugate linear isometric operator $\widehat J$ on $\Cal H$ for the dual $(\widehat M,\widehat\co)$ satisfying $\widehat J\widehat M\widehat J=\widehat M'$ and the modular automorphisms $(\widehat \sigma_t)$ and $({\widehat \sigma}'_t)$ of $\widehat M$, as well as the modular operator $\hat\delta$ for the dual.
\nl
The {\it scaling group} can be characterized as follows:

\inspr{4.7} Proposition \rm
There exists a one-parameter group of automorpisms $(\tau_t)$ of $(M,\co)$ such that 
$$\align \co(\sigma_t(x))&=(\tau_t\ot\sigma_t)\co(x) \\ 
         \co(\sigma'_t(x))&=(\sigma'_t\ot \tau_{-t})\co(x) 
\endalign$$ 
for all $x\in M$ and $t\in \Bbb R$. All the automorphisms in $(\sigma_t)$, $(\sigma'_t)$ and $(\tau_t)$ mutually commute.
\einspr

Similarly, we have the scaling group $(\widehat \tau_t)$ on the dual, characterized by similar formulas.
\snl
If we take a proper analytic extension, we see that $\tau_{-i}$ is like the square $S^2$ of the antipode in a $^*$-algebraic quantum group. The first formula replaces $\co(\sigma(a)))=(S^2\ot\sigma)\co(a)$ and the second one is $\co(\sigma'(a))=(\sigma'\ot S^{-2})\co(a)$ for an element $a$ in a $^*$algebraic quantum group.
\snl
Again, the proof is technically rather difficult. It essentially uses the polar decomposition of an operator $\Lambda(x)\mapsto \Lambda(S(x)^*)$ where $S$ is the 'antipode', roughly defined by the formula
$$S((\iota\ot\varphi)(\co(x)(1\ot y)))=(\iota\ot\varphi)((1\ot x)\co(y))$$
for well-chosen elements $x$ and $y$ in the von Neumann algebra $M$.
\nl
There are several relations among the data we have so far:

\inspr{4.8} Proposition \rm
When $x\in M$ and $y\in \widehat M$ we have
$$\alignat{2} 
\sigma_t(x)&=\nabla^{it}x\nabla^{-it} &\qquad\qquad \tau_t(x)&=\widehat\nabla^{it}x\widehat\nabla^{-it}\\
\widehat\sigma_t(y)&=\widehat\nabla^{it}y\widehat\nabla^{-it} &\qquad\qquad \widehat\tau_t(y)&=\nabla^{it}y\nabla^{-it}
\endalignat$$
for all $t\in \Bbb R$.
\einspr

The formulas on the left were mentioned already but the others are new (and somewhat remarkable). We do not have any counterparts of these equations in the theory of $^*$-algebraic quantum groups. This is not so with the following results.

\inspr{4.9} Proposition \rm
There exists a strictly positive number $\nu$, called the {\it scaling constant}, satisfying
$$\alignat{2}
\varphi\circ\tau_t&=\nu^{-t}\varphi &\qquad\qquad \varphi\circ\sigma'_t&=\nu^t \varphi\\
\psi\circ\tau_t&=\nu^{-t}\psi &\qquad\qquad \psi\circ\sigma_t&=\nu^{-t} \psi
\endalignat$$
for all $t\in \Bbb R$.
\einspr

When extending these formulas analytically to the complex number $-i$, we find e.g.\ $\varphi\circ S^2=\nu^{i}\varphi$ and we see that $\nu^i$ turns out to replace the scaling constant as introduced for $^*$-algebraic quantum groups. As mentioned already, in this case, the scaling constant can be non-trivial, see e.g.\ [VD5].
\snl
Also the above result is a consequence of the uniqueness of the invariant weights.
\snl
And finally, we have some formulas relating $\delta$ with the other data:

\iinspr{4.10} Proposition \rm
We have $\tau_t(\delta)=\delta$ and
$$\sigma_t(\delta)=\nu^t\delta\qquad\qquad \sigma'_t(\delta)=\nu^{-t}\delta$$
for all $t$. We also have $\widehat J \delta \widehat J=\delta^{-1}$.
\einspr

Of course, these formulas have to be interpreted (e.g.\ by looking at powers $\delta^{is}$ of $\delta$). There is also a formula for $J\delta J$ but that is more complicated. Similar equations hold for the dual modular operator $\widehat \delta$.

\nl
Having defined the main objects and the most important formulas, we can now state the analytical form of Radford's formula for locally compact quantum groups (see Theorem 4.20 in [VD8]):

\iinspr{4.11} Theorem \rm Because the left Haar weight is relatively invariant, we can define a one-parameter group of unitary operators, denoted $P^{it}$, by the formula $P^{it}\Lambda(x)=\nu^{\frac12 t}\Lambda(\tau_t(x))$ for all $x\in \Cal N$. Then we have
$$P^{-2it}=\delta^{it}(J\delta^{it}J)\widehat\delta^{it}(\widehat J\widehat\delta^{it}\widehat J)$$
for all $t$.
\einspr

Compare this formula, call it the 'second' formula in what follows, with the formula in Theorem 3.8, which we will call the 'first' one. And assume for the moment that the scaling constant is $1$. Change $t$ to $-t$ in the first formula and 'apply' $\Lambda$. On the left hand side, we get $P^{-2it}\Lambda(a)$. When we look at the right hand side, first we have left multiplication with $\delta^{it}$ in the first formula which we find as the operator $\delta^{it}$ in the second formula. Next we have right multiplication with $\delta^{-it}$ in the first formula that results in the operator $J\delta^{it}J$ in the second formula. The change in sign comes from the fact that $J$ is conjugate linear and $\delta$ self-adjoint. Also remember that $J$ is the unitary part in the polar decomposition of the map $\Lambda(x)\mapsto \Lambda(x^*)$ and the fact that the involution changes the order allows to express right multiplication with elements as operators, using this map. The third factor in the second formula comes from the left action of $\widehat\delta^{-it}$. Now, the difference in sign is coming from the difference in conventions about the dual coproduct. Flipping this coproduct causes $\widehat\delta$ to be replaced by $\widehat\delta^{-1}$. Finally, the right action of $\widehat\delta^{it}$ corresponds with the factor $\widehat J\widehat\delta^{it}\widehat J$. We have the same sign here because it is changed two times for reasons explained earlier. 
\snl
If the scaling constant is not equal to $1$, we get an extra factor on the left because $P^{-2it}\Lambda(x)=\nu^{-t}\Lambda(\tau_t(x))$. This factor will also occur on the right hand side because right multiplication with $\delta^{-it}$ is not exactly the same as $J\delta^{it}J$. There is a factor $\nu^{\frac12 t}$ coming from the commutation rules between the modular operator $\nabla$ and $\delta$ (as $\sigma(\delta)=\frac{1}{\nu}\delta$ in the case of algebraic quantum groups). Similarly, this scalar will pop up when comparing the right action of $\widehat\delta^{it}$  with the factor $\widehat J\widehat\delta^{it}\widehat J$.
\snl
So, we see that the two formulas are completely in accordance with each other and that it is justified to call the formula in Theorem 4.11 above the analytical form of Radford's formula for locally compact quantum groups.
\snl
Also here, it is interesting to look at some special cases. If e.g.\ $\widehat\delta=1$, also in this case we have $\sigma_t=\tau_t$ and $\sigma'_t=\tau_{-t}$ (as for discrete quantum groups). If both modular operators are $1$, then necessarily the scaling group and the modular automorphisms are trival, causing the Haar weights to be traces and $S^2=\iota$.
\nl\nl
%
%
\bf References \rm
\nl

{\bf [A]} E.\ Abe: \it Hopf algebras. \rm Cambridge University Press (1977).
\snl
{\bf [B-B-T]} M.\ Beattie, D.\ Bulacu \& B.\ Torrecillas: {\it Radford's $S^4$ formula for co-Frobenius Hopf algebras}. J.\ Algebra (2006), to appear. 
\snl
{\bf [DC-VD]} K.\ De Commer \& A.\ Van Daele: {\it Multiplier Hopf algebras imbedded in C$^*$-algebraic quantum groups}. Preprint K.U.\ Leuven (2006). Arxiv math.OA/0611872.
\snl
{\bf[D-VD]} L.\ Delvaux \& A. Van Daele: {\it Algebraic quantum hypergroups}. Preprint University of Hasselt and K.U.\ Leuven (2006). math.RA/0606466
\snl
{\bf [D-VD-W]} L.\ Delvaux, A. Van Daele \& S.\ Wang: {\it A note on Radford's $S^4$ formula}. Preprint  University of Hasselt, K.U.\ Leuven and Nanjing University (2006). \newline Arxiv math.RA/0608096.
\snl
{\bf [Dr-VD]} B.\ Drabant \& A. Van Daele: {\it Pairing and Quantum double of multiplier Hopf algebras}.  Algebras and Representation Theory 4 (2001), 109-132.
\snl
{\bf [Dr-VD-Z]} B.\ Drabant, A.\ Van Daele \& Y.\ Zhang: {\it Actions of multiplier Hopf algebras}. Commun.\ Alg.\ 27 (1999), 4117--4172.
\snl
{\bf [E-R]} E.G.\ Effros \& Z.-J.\ Ruan : {\it Discrete quantum groups I.
The Haar measure}. Int.\ J.\ Math.\ 5 (1994) 681-723.
\snl
{\bf [E-S]} M.\ Enock \& J.-M.\ Schwartz: {\it Kac algebras and duality
for locally compact groups.} Springer (1992).
\snl
{\bf [K]} J.\ Kustermans: {\it The analytic structure of algebraic quantum groups}. J.\ Algebra 259 (2003), 415--450.
\snl
{\bf [K-R]} L.\ Kaufman \& D.\ Radford: {\it A necessary and sufficient condition for a finite-dimensional Drinfel'd double to be a ribbon Hopf algebra}. J.\ of Alg.\ 159 (1993), 98-114.
\snl
{\bf [K-V1]} J.\ Kustermans \& S.\ Vaes: {\it A simple definition for locally compact quantum groups}. C.R.\ Acad.\ Sci.\ Paris S\'er I 328 (1999), 871--876.
\snl
{\bf [K-V2]} J.\ Kustermans \& S.\ Vaes: {\it Locally compact quantum groups}. Ann.\ Sci.\ \'Ecole Norm.\ Sup. (4) (33) (2000), 837--934.
\snl
{\bf [K-V3]} J.\ Kustermans \& S.\ Vaes: {\it Locally compact quantum groups in the von Neumann algebra setting}. Math.\ Scand.\ 92 (2003), 68--92. 
\snl
{\bf [K-VD]} J.\ Kustermans \& A.\ Van Daele: {\it C$^*$-algebraic quantum groups arising from algebraic quantum groups}. Int.\ J.\ Math.\ 8 (1997), 1067--1139.
\snl
{\bf [L-VD]} M.B.\ Landstad \& A.\ Van Daele: {\it Groups with compact open subgroups and multiplier Hopf $^*$-algebras}. Preprint University of Trondheim \& K.U.\ Leuven. \newline 
Arxiv math.0A/0702458 (to appear in Expositiones Mathematica).
\snl
{\bf [M-VD]} A.\ Maes \& A.\ Van Daele:  {\it Notes on compact
quantum groups}. Nieuw Archief voor Wiskunde, Vierde serie
16 (1998), 73--112.
\snl
{\bf [M-N]} M.\ Masuda \& Y.\ Nakagami: {\it A von Neumann algebra
framework for the duality of quantum groups.} Publ. RIMS Kyoto
{\bf 30} (1994), 799--850.
\snl
{\bf [M-N-W]} T.\ Masuda, Y.\ Nakagami \& S.L.\ Woronowicz: {\it A C$^*$-algebraic framework for the quantum groups}. Int.\ J.\ of Math.\ 14 (2003), 903--1001.  
\snl
{\bf [P-W]} P.\ Podle\'s \& S.L.\ Woronowicz: {\it Quantum deformation of Lorentz group}. Commun.\ Math.\ Phys.\ 130 (1990), 381--431.
\snl
{\bf [R]} D.\ Radford: {\it The order of the antipode of any finite-dimensional Hopf algebra is finite}. Amer.\ J.\ Math.\ 98 (1976), 333--355. 
\snl
{\bf [S]} M.\ Sweedler: {\it Hopf algebras}. Benjamin, New-York (1969). 
\snl
{\bf [T1]} M.\ Takesaki: {\it Theory of Operator Algebras I}. Springer-Verlag, New York (1979). 
\snl
{\bf [T2]} M.\ Takesaki: {\it Theory of Operator Algebras II}. Springer-Verlag, New York (2001).
\snl
{\bf [VD1]} A.\ Van Daele: {\it Multiplier Hopf algebras}. Trans.\ Am.\ Math.\ Soc.\ 342 (1994), 917--932.
\snl
{\bf [VD2]} A.\ Van Daele: {\it Discrete quantum groups.} J.\ of Alg.\ 180 (1996), 431-444.
\snl
{\bf [VD3]} A.\ Van Daele:  {\it An algebraic framework for group duality}.  Adv.\ Math.\ 140 (1998), 323--366.
\snl
{\bf [VD4]} A.\ Van Daele: {\it Quantum groups with invariant integrals}. Proc.\ Natl.\ Acad.\ Sci.\ USA  97 (2000), 541-546.
\snl
{\bf [VD5]} A.\ Van Daele: {\it The Haar measure on some locally compact quantum groups}. Pre\-print K.U.\ Leuven (2001). Arxiv math.OA/0109004.
\snl
{\bf [VD6]} A.\ Van Daele: {\it Multiplier Hopf $^\ast$-algebras with positive integrals: A laboratory for locally compact quantum groups.} Irma Lectures in Mathematical and Theoretical Physics 2: Locally compact Quantum Groups and Groupoids. Proceedings of the meeting in Strasbourg on Hopf algebras, quantum groups and their applications (2002). Ed.\ V.\ Turaev \& L.\ Vainerman. Walter de Gruyter, (2003), 229--247.
\snl
{\bf [VD7]} A.\ Van Daele: {\it Locally compact quantum groups: The von Neumann algebra versus the C$^*$-algebra approach}. Bulletin of Kerala Mathematics Association, Special issue (2006), 153--177.
\snl
{\bf [VD8]} A.\ Van Daele: {\it Locally compact quantum groups. A von Neumann algebra approach}. Preprint K.U.\ Leuven (2006). Arxiv math.OA/0602212.
\snl
{\bf [VD-Z]} A.\ Van Daele \& Y.\ Zhang: {\it A survey on multiplier Hopf algebras}. In 'Hopf algebras and Quantum Groups', eds. S.\ Caenepeel \& F.\ Van Oyestayen, Dekker, New York (1998), pp. 259--309.
\snl
{\bf [W1]} S.L.\ Woronowicz: {\it Twisted $SU(2)$ group. An example of a non-commutative differential calculus}. Publ.\ RIMS, Kyoto University 23 (1987), 117--181.
\snl
{\bf [W2]} S.L.\ Woronowicz: {\it Compact matrix pseudogroups}. Comm. Math. Phys. 111 (1987), 613-665.
\snl
{\bf [W3]} S.L.\ Woronowicz: {\it  Compact quantum groups.}
Quantum symmetries/Symm\'{e}tries quantiques.  Proceedings of the
Les Houches summer school 1995, North-Holland, Amsterdam (1998),
845--884.
\snl

\end